\documentclass[12pt]{article}
\usepackage{amssymb}
\usepackage{amsfonts}
\usepackage{times}
\usepackage{mathptmx}
\usepackage{amsmath}
\usepackage[usenames]{color}
\usepackage{mathrsfs}
\usepackage{amsfonts}
\usepackage{amssymb,amsmath}
\usepackage{CJK}
\usepackage{cite}
\usepackage{cases}
\usepackage{amsthm}
\def\cl{\centerline}

\def\vs{\vspace*}

\def\QED{\hfill$\Box$}
\def\ni{\noindent}

\numberwithin{equation}{section}
\newtheorem{theo}{Theorem}[section]
\newtheorem{defi}[theo]{Definition}

\newtheorem{prop}[theo]{Proposition}

\begin{document}
\begin{center}
{\bf\large Generalized derivations of n-Hom Lie superalgebras}
\footnote {Supported by the National Natural Science Foundation of China (No.~11431010).

$^{\,\S}$Corresponding author: G.~Fan.
}
\end{center}

\cl{Jinsen Zhou$^{\,*}$, Guangzhe Fan$^{\,\S}$
}

\cl{\small $^{\,*}$zjs9932@126.com}
\cl{\small $^{\,*}$School of Information Engineering, Longyan University,
 Longyan 364012, Fujian, P. R. China}
\cl{\small $^{\,\S}$yzfanguangzhe@126.com}
\cl{\small $^{\,\S}$Department of Mathematics, Tongji University, Shanghai 200092,
P. R. China}

\vs{8pt}

{\small\footnotesize
\parskip .005 truein
\baselineskip 3pt \lineskip 3pt
\noindent{{\bf Abstract:} It is well known that n-Hom Lie superalgebras are certain generalizations of n-Lie algebras. This paper is devoted to investigate the generalized derivations of multiplicative n-Hom Lie superalgebras. We generalize the main results of Leger and Luks to the case of multiplicative n-Hom Lie superalgebras. Firstly, we review some concepts associated with a multiplicative n-Hom Lie superalgebra $N$. Furthermore, we give the definitions of the generalized derivations, quasiderivations, center derivations, centroids and quasicentroids. Obviously, we have the following tower
$ZDer(N)\subseteq Der(N)\subseteq QDer(N)\subseteq GDer(N)\subseteq End(N)$.
Later on, we give some useful properties and connections between these derivations. Moreover, we obtain that the quasiderivation of $N$ can be embedded as a derivation in a larger multiplicative n-Hom Lie superalgebra. Finally, we conclude that the derivation of the larger multiplicative n-Hom Lie superalgebra has a direct sum decomposition when the center of $N$ is equal to zero.
\vs{5pt}


}}
\parskip .001 truein\baselineskip 6pt \lineskip 6pt

\section{Introduction}

In 1973, the first instance of n-ary algebra was proposed by Nambu \cite{N} as a generalization of the Hamiltonian mechanics. In 1985, in \cite{F} Filippov  introduced the abstract definitions of n-algebras or n-Lie algebras (when the bracket is skew-symmetric). He also classified n-Lie algebras of $n+1$ dimension over an algebraically closed field of characteristic zero. In \cite{K}, Kasymov investigated the structure and representation theories of n-Lie algebras.
Recently, the structure theories and representation theories related to n-Lie algebras have been extensively investigated in \cite{KP,AMM,BL}. In this paper, we introduce the concept of  n-Hom Lie superalgebras as generalizations of n-Lie algebras. We shall investigate the generalized derivations of multiplicative n-Hom Lie superalgebras, which would lead to promote the development of structure theories of n-Lie (super)algebras.

As is well known, derivations and generalized derivations are very important objects in the research of
 Lie algebras and its generalizations. The most important and systematic research on the generalized derivations
 of a Lie algebra and their subalgebras was due to Leger and Luks \cite{LL}. In their articles, they obtained many nice properties of generalized derivations.
 Their results were generalized by many authors. For example, Zhang and Zhang in \cite{ZZ} generalized the above results to the case of Lie superalgebras; Chen, Zhou, etc considered the generalized derivations of color Lie algebras, Hom-Lie (super)algebras, Lie triple systems and Hom-Lie triple systems \cite{CMN,ZNC,ZCM1,ZCM2,ZCM3}. Later on, in \cite{KP} Kaygorodov and Popov generalized the results of Leger and Luks about generalized derivations of Lie algebras to the case of (color) n-ary algebras. In a sense, the multiplicative n-Hom Lie superalgebras are generalizations of those algebras. Hence, some results about generalized derivations of those algebras can be applied in this case.

 The present paper is devoted to study all kinds of generalized derivations of multiplicative n-Hom Lie superalgebras. We mainly investigate the derivation algebra
$Der(N)$, the center derivation algebra $ZDer(N)$, the generalized derivation algebra $GDer(N)$, and the quasiderivation algebra $QDer(N)$ of a multiplicative n-Hom Lie superalgebra $N$. Furthermore, we give some useful properties and connections between these derivations. 

Here is a detailed outline of the contents of the main parts of the article. In Section 2, we recall some basic definitions and notations of multiplicative n-Hom Lie superalgebras used in this paper. In Section 3, we give some elementary observations about generalized derivations, quasiderivations, centroids and quasicentroids, some of which are technical results to be in the sequel. In Section 4, we prove that the quasiderivation of $N$ can be embedded as a derivation in a larger multiplicative n-Hom Lie superalgebra. Finally, we obtain that $Der(\breve{N})$ has a direct sum decomposition when the center of $N$ is equal to zero.

Throughout this paper, we denote by $\mathbb{F}$, $\mathbb{Z}_{2}$, $\mathbb{N}$ a field of characteristic zero, the ring of integers modulo 2, the set of non-negative integers. Set $J=\{1,2,\cdots,n\}$.

\section{Preliminaries}

In this section we review the theory of n-Hom Lie superalgebras and generalize some results of \cite{ZCM1,KP}. We give the definitions of generalized derivations which we refer to \cite{LL}. For readers' convenience, we introduce some notations used in this paper.

Let $V$ be a vector superspace over $\mathbb{F}$ that is a  $\mathbb{Z}_{2}$-graded vector space with a direct sum $V=V_{\overline{0}}\oplus V_{\overline{1}}$. An elements $x\in V_{\gamma}\  (\gamma \in \mathbb{Z}_{2})$ is said to be homogeneous of degree $\gamma$. Denote $hg(V)$ by the homogeneous elements of $V$.
The degree of an element $x$ is denoted by $|x|$. Let $V$ and
$W$ be two $\mathbb{Z}_{2}$-graded vector spaces.
The vector space of all linear maps from $V$ to $W$ is denoted by $Hom(V,W)=Hom(V,W)_{\overline{0}}\oplus Hom(V,W)_{\overline{1}}$. A linear map $f:V\rightarrow W$ is said to be homogeneous of
 degree $\xi\in \mathbb{Z}_{2}$, if $f(x)$ is homogeneous of degree $\gamma+\xi$ for all the element
$x\in V_{\gamma}$. The set of all such maps is denoted by $Hom(V,W)_{\xi}$. In addition, we call $f$ an even linear map  from $V$ to $W$, if $f(V_{\gamma})\subseteq W_{\gamma}$ holds for any $\gamma \in \mathbb{Z}_{2}$. The vector space of all linear maps from $V$ to $V$ is denoted by $End(V)$. The notations $\xi,\eta,\theta,\gamma$ denote the elements of $\mathbb{Z}_{2}$ unless otherwise stated. Moreover, $x$ is always assumed to be homogeneous when $|x|$ occurs. Set $|X_{i}|=|x_{1}|+\cdots+|x_{i}|$. In particular, we set $|X_{0}|=0$.

\begin{defi}\label{250}\rm
An \emph{n-Lie superalgebra} is a pair $(N,[\cdot,\cdots,\cdot])$ consisting of a $\mathbb{Z}_{2}$-graded vector space $N=N_{\overline{0}}\oplus N_{\overline{1}}$ and a multilinear map
$[\cdot,\cdots,\cdot]: \underbrace{N\times N \times \cdots \times N}_{n}\rightarrow N$, satisfying
\begin{equation*}\label{251}
|[x_{1},\cdots,x_{n}]|=|X_{n}|,
\end{equation*}
\begin{equation*}\label{252}
[x_{1},\cdots,x_{i},x_{i+1},\cdots,x_{n}]=-(-1)^{|x_{i}||x_{i+1}|}[x_{1},\cdots,x_{i+1},x_{i},\cdots,x_{n}],
\end{equation*}
\begin{equation*}\label{253}
[x_{1},\cdots,x_{n-1},[y_{1},\cdots,y_{n}]]=\sum^{n}_{i=1}(-1)^{|X_{n-1}||Y_{i-1}|}
[y_{1},\cdots,y_{i-1},[x_{1},\cdots,x_{n-1},y_{i}],y_{i+1},\cdots,y_{n}],
\end{equation*}
for any $x_{i},y_{i} \in hg(N), i\in J$.

When $N_{\overline{1}}=\{0\}$, then $N$ is actually an n-Lie algebra.
\end{defi}

\begin{defi}\rm
An \emph{n-Hom Lie superalgebra} is $(N,[\cdot,\cdots,\cdot],\alpha_{1},\cdots,\alpha_{n-1})$ consisting of a $\mathbb{Z}_{2}$- graded vector space $N=N_{\overline{0}}\oplus N_{\overline{1}}$ and a multilinear map
$[\cdot,\cdots,\cdot]:\underbrace{N\times N \times \cdots \times N}_{n}\rightarrow N$ and a family $(\alpha_{i})_{1\leq i\leq n-1}$ of even linear maps $\alpha_{i}:N\rightarrow N$, satisfying
\begin{equation*}\label{251}
|[x_{1},\cdots,x_{n}]|=|X_{n}|,
\end{equation*}
\begin{equation*}\label{252}
[x_{1},\cdots,x_{i},x_{i+1},\cdots,x_{n}]=-(-1)^{|x_{i}||x_{i+1}|}[x_{1},\cdots,x_{i+1},x_{i},\cdots,x_{n}],
\end{equation*}
\begin{equation*}\label{254}
[\alpha_{1}(x_{1}),\cdots,\alpha_{n-1}(x_{n-1}),[y_{1},\cdots,y_{n}]]
\end{equation*}
\begin{equation*}\label{255}
=\sum^{n}_{i=1}(-1)^{|X_{n-1}||Y_{i-1}|}
[\alpha_{1}(y_{1}),\cdots,\alpha_{i-1}(y_{i-1}),[x_{1},\cdots,x_{n-1},y_{i}],
\alpha_{i}(y_{i+1}),\cdots,\alpha_{n-1}(y_{n})],
\end{equation*}
for any $x_{i},y_{i} \in hg(N), i\in J$.
\end{defi}
We also see that if $\alpha_{1}=\cdots=\alpha_{n-1}=id_{N}$, then $N$ is just an n-Lie superalgebra.

\begin{defi}\rm
An \emph{n-Hom Lie superalgebra} $(N,[\cdot,\cdots,\cdot],\alpha_{1},\cdots,\alpha_{n-1})$ is \emph{multiplicative}, if
$(\alpha_{i})_{1\leq i\leq n-1}$ with $\alpha_{1}=\cdots=\alpha_{n-1}=\alpha$ and satisfying

$$\alpha([x_{1},\cdots,x_{n}])=[\alpha(x_{1}),\cdots,\alpha(x_{n})],$$
\begin{equation*}\label{254}
[\alpha(x_{1}),\cdots,\alpha(x_{n-1}),[y_{1},\cdots,y_{n}]]
\end{equation*}
\begin{equation*}\label{256}
=\sum^{n}_{i=1}(-1)^{|X_{n-1}||Y_{i-1}|}
[\alpha(y_{1}),\cdots,\alpha(y_{i-1}),[x_{1},\cdots,x_{n-1},y_{i}],
\alpha(y_{i+1}),\cdots,\alpha(y_{n})],
\end{equation*}
for any $x_{i},y_{i} \in hg(N), i\in J$.
\end{defi}

For convenience, from now on, we always assume that $(N,[\cdot,\cdots,\cdot],\alpha)$ is a multiplicative n-Hom Lie superalgebra over $\mathbb{F}$ unless otherwise stated.

\begin{defi}\rm
Define the
 following vector subspace $\Omega$ of End($N$) consisting of linear maps on $N$ as following:
 \begin{eqnarray*}\!\!\!\!\!\!\!\!\!\!\!\!&\!\!\!\!\!\!\!\!\!\!\!\!\!\!\!&
\ \ \ \ \ \ \ \ \Omega=\{u\in End(N)\mid u \alpha=\alpha u\}.
\end{eqnarray*}
 and
\begin{eqnarray*}\!\!\!\!\!\!\!\!\!\!\!\!&\!\!\!\!\!\!\!\!\!\!\!\!\!\!\!&
\ \ \ \ \ \ \ \ \widetilde{\alpha}:\Omega\rightarrow \Omega;\ \ \widetilde{\alpha}(u)=\alpha u.
\end{eqnarray*}
Then $(\Omega,[\cdot,\cdot],\widetilde{\alpha})$ is a Hom-Lie superalgebra over $\mathbb{F}$ and satisfies the following bracket:
\begin{eqnarray*}\!\!\!\!\!\!\!\!\!\!\!\!&\!\!\!\!\!\!\!\!\!\!\!\!\!\!\!&
\ \ \ \ \ \ \ \ [D_{\xi},D_{\eta}]=D_{\xi}D_{\eta}-(-1)^{\xi\eta}D_{\eta}D_{\xi},
\end{eqnarray*}
for any $D_{\xi},D_{\eta}\in hg(\Omega)$.
 %%All homogeneous elements of $\Omega$ are denoted by $hg(\Omega)$.
\end{defi}

\begin{defi}\rm
Let $(L,[\cdot,\cdot],\alpha)$ be a Hom-Lie superalgebra.
 A graded subspace $M\subseteq L$ is a Hom-subalgebra of $L$ if $\alpha(M)\subseteq M$ and $M$ is closed under the bracket operation $[\cdot,\cdot]$, i.e. $[D_{\xi},D_{\eta}]\in M$, for any $D_{\xi},D_{\eta}\in M$.

A graded subspace $I\subseteq L$ is called a Hom-ideal of $(L,[\cdot,\cdot],\alpha)$ if $\alpha(I)\subseteq I$ and $I$ is closed under the bracket operation $[\cdot,\cdot]$, i.e. $[D_{\xi},D_{\eta}]\in I$, for any $D_{\xi}\in L, D_{\eta}\in I$.
\end{defi}

\begin{defi}\rm
A homogeneous
linear map $D:N\rightarrow N$ of degree $\xi$ is said to be an $\alpha^{k}$-derivation of $N$ for $k\in \mathbb{N}$, if it satisfies
$$D\alpha=\alpha D,$$
$$D([x_{1},\cdots,x_{n}])=\sum^{n}_{i=1}(-1)^{\xi |X_{i-1}|}
[\alpha^{k}(x_{1}),\cdots,\alpha^{k}(x_{i-1}),D(x_{i}),\alpha^{k}(x_{i+1}),\cdots,\alpha^{k}(x_{n})],$$
for any $x_{i}\in hg(N),i\in J$.
\end{defi}

We denote by $Der_{\alpha^{k}}(N)$ the set of all $\alpha^{k}$-derivations, then $Der(N):=\oplus_{k\geq 0}Der_{\alpha^{k}}(N)$ provided with the super-commutator and the following even map
$$\widetilde{\alpha}:Der(N)\rightarrow Der(N);~~~~\ \ \widetilde{\alpha}(D)=D\alpha$$
is a Hom-subalgebra of $\Omega$ and is called the derivation algebra of $N$.~(see \cite{GC})

\begin{defi}\rm
An endomorphism
$D\in End_{\xi}(N)$ is said to be a homogeneous generalized $\alpha^{k}$-derivation of degree $\xi$ of $N$, if
 there exist endomorphisms $D^{(i)}\in End_{\xi}(N)$ such that
 $$D\alpha=\alpha D,D^{(i)}\alpha=\alpha D^{(i)},$$
 $$[D(x_{1}),\alpha^{k}(x_{2}),\cdots,\alpha^{k}(x_{n})]+\sum^{n}_{i=2}(-1)^{\xi |X_{i-1}|}
 [\alpha^{k}(x_{1}),\cdots,\alpha^{k}(x_{i-1}),D^{(i-1)}(x_{i}),\alpha^{k}(x_{i+1}),\cdots,\alpha^{k}(x_{n})]$$
 $$=D^{(n)}([x_{1},\cdots,x_{n}]),$$
for any $x_{i}\in hg(N),i\in J$.
\end{defi}

\begin{defi}\rm
An endomorphism
$D\in End_{\xi}(N)$ is called a homogeneous $\alpha^{k}$-quasiderivation of degree $\xi$ of $N$, if
 there exists $D^{\prime}\in End_{\xi}(N)$ such that
$$D\alpha=\alpha D,D^{\prime}\alpha=\alpha D^{\prime},$$
$$[D(x_{1}),\alpha^{k}(x_{2}),\cdots,\alpha^{k}(x_{n})]+\sum^{n}_{i=2}(-1)^{\xi |X_{i-1}|}
 [\alpha^{k}(x_{1}),\cdots,\alpha^{k}(x_{i-1}),D(x_{i}),\alpha^{k}(x_{i+1}),\cdots,\alpha^{k}(x_{n})]$$
 $$=D^{\prime}([x_{1},\cdots,x_{n}]),$$
for any $x_{i}\in hg(N),i\in J$.
\end{defi}

Let $GDer_{\alpha^{k}}(N)$ and $QDer_{\alpha^{k}}(N)$ be the set of homogeneous generalized $\alpha^{k}$-derivations and of homogeneous $\alpha^{k}$-quasiderivation, respectively. That is
$$GDer(N):=\oplus_{k\geq 0}GDer_{\alpha^{k}}(N),\ \ QDer(N):=\oplus_{k\geq 0}QDer_{\alpha^{k}}(N).$$

\begin{defi}\rm
A homogeneous
linear map $D:N\rightarrow N$ of degree $\xi$ is said to be an $\alpha^{k}$-centroid of $N$ for $k\in \mathbb{N}$, if it satisfies
the following equations:
$$D\alpha=\alpha D,$$
$$[D(x_{1}),\alpha^{k}(x_{2}),\cdots,\alpha^{k}(x_{n})]=(-1)^{\xi |X_{i-1}|}
 [\alpha^{k}(x_{1}),\cdots,\alpha^{k}(x_{i-1}),D(x_{i}),\alpha^{k}(x_{i+1}),\cdots,\alpha^{k}(x_{n})]$$
 $$=D([x_{1},\cdots,x_{n}]),$$
for any $x_{i}\in hg(N),i\in J$.
\end{defi}

\begin{defi}\rm
We call $D\in End_{\xi}(N)$  a homogeneous $\alpha^{k}$-quasicentroid of degree $\xi$ of $N$, if
it satisfies
$$D\alpha=\alpha D,$$
$$[D(x_{1}),\alpha^{k}(x_{2}),\cdots,\alpha^{k}(x_{n})]=(-1)^{\xi |X_{i-1}|}
 [\alpha^{k}(x_{1}),\cdots,\alpha^{k}(x_{i-1}),D(x_{i}),\alpha^{k}(x_{i+1}),\cdots,\alpha^{k}(x_{n})],$$
for any $x_{i}\in hg(N),i\in J$.
\end{defi}

Denote by $C_{\alpha^{k}}(N)$, $QC_{\alpha^{k}}(N)$ the set of homogeneous $\alpha^{k}$-centroid, $\alpha^{k}$-quasicentroid, respectively. That is,
$$C(N):=\oplus_{k\geq 0}C_{\alpha^{k}}(N),\ \ QC(N):=\oplus_{k\geq 0}QC_{\alpha^{k}}(N).$$

\begin{defi}\rm
 An endomorphism
$D\in End_{\xi}(N)$ is said to be a homogeneous $\alpha^{k}$-center derivations of degree $\xi$ of $N$, if
 it satisfies
$$D\alpha=\alpha D,$$
$$[D(x_{1}),\alpha^{k}(x_{2}),\cdots,\alpha^{k}(x_{n})]=D([x_{1},\cdots,x_{n}])=0,$$
for any $x_{i}\in hg(N),i\in J$.
\end{defi}

Denote by $ZDer_{\alpha^{k}}(N)$ the set of homogeneous $\alpha^{k}$-center derivations. That is
$$ZDer(N):=\oplus_{k\geq 0}ZDer_{\alpha^{k}}(N).$$

Obviously, we have the following tower
$$ZDer(N)\subseteq Der(N)\subseteq QDer(N)\subseteq GDer(N)\subseteq End(N).$$

\begin{defi}\rm
 If $Z(N):=Z_{\overline{0}}(N)\oplus Z_{\overline{1}}(N)$, with
 $Z_{\xi}(N)=\{x\in hg_{\xi}(N)\ |\ [x,y_{2},y_{3},\cdots,y_{n}]=0, \forall\  y_{2},y_{3},\cdots,y_{n}\in N\},$ then $Z(N)$ is called the center of $N$.
\end{defi}

\section{Generalized derivation algebras and their Hom-subalgebras}

In this section, we investigate some basic properties of generalized derivations, quasiderivations and center derivations of a multiplicative n-Hom Lie superalgebra. Let $(N,[\cdot,\cdots,\cdot],\alpha)$ is a multiplicative n-Hom Lie superalgebra over $\mathbb{F}$ unless otherwise stated.

\begin{prop}\label{po1}
We get the following statements:

(1) $GDer(N), QDer(N)$ and $C(N)$ are Hom-subalgebras of $\Omega$;

(2) $ZDer(N)$ is a Hom-ideal of $Der(N)$.
\end{prop}
\ni\ni{\it Proof.}\ \
(1) Assume that $D_{\xi}\in GDer_{\alpha^{k}}(N),\ D_{\eta}\in GDer_{\alpha^{s}}(N)$ for any $x_{i} \in hg(N)$, where $k,s\in\mathbb{N}$, $i\in J$, then we have

$[\widetilde{\alpha}(D_{\xi})(x_{1}),\alpha^{k+1}(x_{2}),\cdots,\alpha^{k+1}(x_{n})]$

$=[(D_{\xi}\alpha)(x_{1}),\alpha^{k+1}(x_{2}),\cdots,\alpha^{k+1}(x_{n})]$

$=\alpha([D_{\xi}(x_{1}),\alpha^{k}(x_{2}),\cdots,\alpha^{k}(x_{n})])$

$=\alpha(D_{\xi}^{(n)}[x_{1},x_{2},\cdots,x_{n}]
-\sum^{n}_{i=2}(-1)^{\xi |X_{i-1}|}
 [\alpha^{k}(x_{1}),\cdots,\alpha^{k}(x_{i-1}),D^{(i-1)}_{\xi}(x_{i}),
 \alpha^{k}(x_{i+1}),\cdots,\alpha^{k}(x_{n})])$

$=\widetilde{\alpha}(D_{\xi}^{(n)})([x_{1},x_{2},\cdots,x_{n}])$

$-\sum^{n}_{i=2}(-1)^{\xi |X_{i-1}|}
 [\alpha^{k+1}(x_{1}),\cdots,\alpha^{k+1}(x_{i-1}),\widetilde{\alpha}(D^{(i-1)}_{\xi})(x_{i}),
 \alpha^{k+1}(x_{i+1}),\cdots,\alpha^{k+1}(x_{n})]$.

For any $i\in J$, we can easily see that $\widetilde{\alpha}(D^{(i)}_{\xi})$ and  $\widetilde{\alpha}(D_{\xi})$ belong to $End_{\xi}(N)$ and
$GDer_{\alpha^{k+1}}(N)$ of degree $\xi$ respectively,
then we also get

$[D_{\xi}D_{\eta}(x_{1}),\alpha^{k+s}(x_{2}),\cdots,\alpha^{k+s}(x_{n})]$

$=D_{\xi}^{(n)}([D_{\eta}(x_{1}),\alpha^{s}(x_{2}),\cdots,\alpha^{s}(x_{n})])$

$-\sum_{i=2}^{n}(-1)^{\xi (\eta+|X_{i-1}|)}[\alpha^{k}(D_{\eta}(x_{1})),\alpha^{k+s}(x_{2}),\cdots,\alpha^{k+s}(x_{i-1}),
D_{\xi}^{(i-1)}(\alpha^{s}(x_{i})),\alpha^{k+s}(x_{i+1}),\cdots,\alpha^{k+s}(x_{n})]$

$=D_{\xi}^{(n)}D_{\eta}^{(n)}([x_{1},x_{2},\cdots,x_{n}])$

$-\sum^{n}_{j=2}(-1)^{\eta |X_{j-1}|}D_{\xi}^{(n)}
[\alpha^{s}(x_{1}),\cdots,\alpha^{s}(x_{j-1}),D_{\eta}^{(j-1)}(x_{j}),
\alpha^{s}(x_{j+1}),\cdots,\alpha^{s}(x_{n})]$

$-\sum^{n}_{i=2}(-1)^{\xi (\eta+|X_{i-1}|)}D_{\eta}^{(n)}
[\alpha^{k}(x_{1}),\cdots,\alpha^{k}(x_{i-1}),D^{(i-1)}_{\xi}(x_{i}),
 \alpha^{k}(x_{i+1}),\cdots,\alpha^{k}(x_{n})]$

$+\sum^{n}_{i=2,j=2,j<i}(-1)^{\xi (\eta+|X_{i-1}|)+\eta |X_{j-1}|}
[\alpha^{k+s}(x_{1}),\cdots,D_{\eta}^{(j-1)}(\alpha^{k}(x_{j})),
\cdots,D_{\xi}^{(i-1)}(\alpha^{s}(x_{i})),\cdots,\alpha^{k+s}(x_{n})]$

$+\sum^{n}_{i=2,j=2,i<j}(-1)^{\xi |X_{i-1}|+\eta |X_{j-1}|}
[\alpha^{k+s}(x_{1}),\cdots,D_{\xi}^{(i-1)}(\alpha^{s}(x_{i})),
\cdots,D_{\eta}^{(j-1)}(\alpha^{k}(x_{j})),\cdots,\alpha^{k+s}(x_{n})]$

$+\sum^{n}_{i=2}(-1)^{\xi (\eta+|X_{i-1}|)+\eta |X_{i-1}|}
[\alpha^{k+s}(x_{1}),\cdots,\alpha^{k+s}(x_{i-1}),
D_{\eta}^{(i-1)}D_{\xi}^{(i-1)}(x_{i}),\alpha^{k+s}(x_{i+1}),\cdots,\alpha^{k+s}(x_{n})].$

Furthermore, one follows that

$[[D_{\xi},D_{\eta}](x_{1}),\alpha^{k+s}(x_{2}),\cdots,\alpha^{k+s}(x_{n})]$

$=[(D_{\xi}D_{\eta}-(-1)^{\xi\eta}D_{\eta}D_{\xi})(x_{1}),\alpha^{k+s}(x_{2}),\cdots,\alpha^{k+s}(x_{n})]$

$=[D_{\xi}^{(n)},D_{\eta}^{(n)}]([x_{1},x_{2},\cdots,x_{n}])
-\sum^{n}_{i=2}(-1)^{(\xi +\eta)|X_{i-1}|}
[\alpha^{k+s}(x_{1}),\cdots,[D_{\xi}^{(i-1)},D_{\eta}^{(i-1)}](x_{i}) \cdots,\alpha^{k+s}(x_{n})].$

For any $2\leq i\leq n$, $[D_{\xi}^{(i-1)},D_{\eta}^{(i-1)}]$ are in $End_{\xi+\eta}(N)$. We know that $[D_{\xi}^{(n)},D_{\eta}^{(n)}]\in GDer_{\alpha^{k+s}}(N)$. Thus, $GDer(N)$ is a Hom-subalgebra of $\Omega$.

Similarly, we also obtain that $QDer(N)$ is a Hom-subalgebra of $\Omega$.

Assume that $D_{\xi}\in C_{\alpha^{k}}(N),\ D_{\eta}\in C_{\alpha^{s}}(N)$ for any  $x_{i} \in hg(N),i\in J$, then

$\widetilde{\alpha}(D_{\xi})([x_{1},\cdots,x_{n}])$

$=\alpha D_{\xi}([x_{1},\cdots,x_{n}])$

$=\alpha([D_{\xi}(x_{1}),\alpha^{k}(x_{2}),\cdots,\alpha^{k}(x_{n})])$

$=[\widetilde{\alpha}(D_{\xi})(x_{1}),\alpha^{k+1}(x_{2}),\cdots,\alpha^{k+1}(x_{n})]$

and

$[\widetilde{\alpha}(D_{\xi})(x_{1}),\alpha^{k+1}(x_{2}),\cdots,\alpha^{k+1}(x_{n})]$

$=\alpha([D_{\xi}(x_{1}),\alpha^{k}(x_{2}),\cdots,\alpha^{k}(x_{n})])$

$=(-1)^{\xi |X_{i-1}|}\alpha
([\alpha^{k}(x_{1}),\cdots,\alpha^{k}(x_{i-1}),D_{\xi}(x_{i}),
\alpha^{k}(x_{i+1}),\cdots,\alpha^{k}(x_{n})])$

$=(-1)^{\xi |X_{i-1}|}
[\alpha^{k+1}(x_{1}),\cdots,\alpha^{k+1}(x_{i-1}),\widetilde{\alpha}(D_{\xi})(x_{i}),
\alpha^{k+1}(x_{i+1}),\cdots,\alpha^{k+1}(x_{n})].$

Hence, we have
\begin{equation*}
\widetilde{\alpha}(D_{\xi})\in C_{\alpha^{k+1}}(N).
\end{equation*}

Note that

$[[D_{\xi},D_{\eta}](x_{1}),\alpha^{k+s}(x_{2}),\cdots,\alpha^{k+s}(x_{n})]$

$=[D_{\xi}D_{\eta}(x_{1}),\alpha^{k+s}(x_{2}),\cdots,\alpha^{k+s}(x_{n})]
-(-1)^{\xi\eta}[D_{\eta}D_{\xi}(x_{1}),\alpha^{k+s}(x_{2}),\cdots,\alpha^{k+s}(x_{n})]$

$=D_{\xi}D_{\eta}([x_{1},x_{2},\cdots,x_{n}])-(-1)^{\xi\eta}D_{\eta}D_{\xi}([x_{1},x_{2},\cdots,x_{n}])$

$=[D_{\xi},D_{\eta}]([x_{1},x_{2},\cdots,x_{n}])$.

Similarly, we have

$(-1)^{(\xi+\eta)|X_{i-1}|}
[\alpha^{k+s}(x_{1}),\cdots,[D_{\xi},D_{\eta}](x_{i}),\cdots,\alpha^{k+s}(x_{n})]
=[D_{\xi},D_{\eta}]([x_{1},x_{2},\cdots,x_{n}]).$

Then, $[D_{\xi},D_{\eta}]\in C_{\alpha^{k+s}}(N)$ of degree $\xi+\eta$. Thus, $C(N)$ is a Hom-subalgebra of $\Omega$.

(2) Assume that $D_{\xi}\in ZDer_{\alpha^{k}}(N),\ D_{\eta}\in Der_{\alpha^{s}}(N)$ for any  $x_{i}\in hg(N), i\in J$, then

$[\widetilde{\alpha}(D_{\xi})(x_{1}),\alpha^{k+1}(x_{2}),\cdots,\alpha^{k+1}(x_{n})]$

$=[(D_{\xi}\alpha)(x_{1}),\alpha^{k+1}(x_{2}),\cdots,\alpha^{k+1}(x_{n})]$

$=\alpha([D_{\xi}(x_{1}),\alpha^{k}(x_{2}),\cdots,\alpha^{k}(x_{n})])$

$=\alpha D_{\xi}([x_{1},x_{2},\cdots,x_{n}])$

$=\widetilde{\alpha}(D_{\xi})([x_{1},x_{2},\cdots,x_{n}])$

$=0$.

Hence, we get
\begin{equation*}
\widetilde{\alpha}(D_{\xi})\in ZDer_{\alpha^{k+1}}(N).
\end{equation*}

Note that

$[D_{\xi},D_{\eta}]([x_{1},x_{2},\cdots,x_{n}])$

$=D_{\xi}D_{\eta}([x_{1},x_{2},\cdots,x_{n}])-(-1)^{\xi\eta}D_{\eta}D_{\xi}([x_{1},x_{2},\cdots,x_{n}])$

$=D_{\xi}([D_{\eta}(x_{1}),\alpha^{s}(x_{2}),\cdots,\alpha^{s}(x_{n})]$

$+\sum^{n}_{i=2}(-1)^{\eta |X_{i-1}|})
[\alpha^{s}(x_{1}),\cdots,\alpha^{s}(x_{i-1}),D^{(i-1)}_{\eta}(x_{i}),
\alpha^{s}(x_{i+1}),\cdots,\alpha^{s}(x_{n})])$

$=0$.

and

$[[D_{\xi},D_{\eta}](x_{1}),\alpha^{k+s}(x_{2}),\cdots,\alpha^{k+s}(x_{n})]$

$=[D_{\xi}D_{\eta}(x_{1}),\alpha^{k+s}(x_{2}),\cdots,\alpha^{k+s}(x_{n})]
-(-1)^{\xi\eta}[D_{\eta}D_{\xi}(x_{1}),\alpha^{k+s}(x_{2}),\cdots,\alpha^{k+s}(x_{n})]$

$=-(-1)^{\xi\eta}[D_{\eta}(D_{\xi}(x_{1})),\alpha^{k+s}(x_{2}),\cdots,\alpha^{k+s}(x_{n})]$

$=-(-1)^{\xi\eta}([D_{\eta}([D_{\xi}(x_{1}),\alpha^{k}(x_{2}),\cdots,\alpha^{k}(x_{n})])$

$-\sum^{n}_{i=2}(-1)^{\eta(\xi+|X_{i-1}|)}
[\alpha^{s}(D_{\xi}(x_{1})),\alpha^{k+s}(x_{2}),\cdots,
D_{\eta}^{(i-1)}(\alpha^{k}(x_{i})),\cdots,\alpha^{k+s}(x_{n})])$

$=0$.

Then, $[D_{\xi},D_{\eta}]\in ZDer_{\alpha^{k+s}}(N)$ of degree $\xi+\eta$. Thus, $ZDer(N)$ is a Hom-ideal of $Der(N)$.   \QED

\begin{prop}\label{po1}
We have the following six statements:

(1) $[Der(N),C(N)]\subseteq C(N)$;

(2) $[QDer(N),QC(N)]\subseteq QC(N)$;

(3) $C(N)\cdot Der (N)\subseteq  Der (N)$;

(4) $C(N)\subseteq  QDer(N)$;

(5) $[QC(T),QC(T)]\subseteq QDer(T)$;

(6) $QDer(N)+QC(N)\subseteq GDer(N)$.
\end{prop}
\ni\ni{\it Proof.}\ \
(1) Suppose that $D_{\xi}\in Der_{\alpha^{k}}(N),\ D_{\eta}\in C_{\alpha^{s}}(N)$ for any $x_{i} \in hg(N), i\in J$.

On one hand, it has

$[D_{\xi}D_{\eta}(x_{1}),\alpha^{k+s}(x_{2}),\cdots,\alpha^{k+s}(x_{n})]$

$=D_{\xi}([D_{\eta}(x_{1}),\alpha^{s}(x_{2}),\cdots,\alpha^{s}(x_{n})])$

$-\sum^{n}_{i=2}(-1)^{\xi(\eta+ |X_{i-1}|)}
[D_{\eta}(\alpha^{k}(x_{1})),\alpha^{k+s}(x_{2}),\cdots,D_{\xi}(\alpha^{s}(x_{i})),\cdots,\alpha^{k+s}(x_{n})]$

$=D_{\xi}D_{\eta}([x_{1},x_{2},\cdots,x_{n}])$.

$-\sum^{n}_{i=2}(-1)^{\xi\eta+(\xi+\eta)|X_{i-1}|}
[\alpha^{k+s}(x_{1}),\cdots,D_{\eta}D_{\xi}(x_{i}),\cdots,\alpha^{k+s}(x_{n})].$

On the other hand, one shows that

$[D_{\eta}D_{\xi}(x_{1}),\alpha^{k+s}(x_{2}),\cdots,\alpha^{k+s}(x_{n})]$

$=D_{\eta}([D_{\xi}(x_{1}),\alpha^{k}(x_{2}),\cdots,\alpha^{k}(x_{n})])$

$=D_{\eta}D_{\xi}([x_{1},x_{2},\cdots,x_{n}])$

$-\sum^{n}_{i=2}(-1)^{\xi|X_{i-1}|}D_{\eta}
([\alpha^{k}(x_{1}),\cdots,D_{\xi}(x_{i}),\cdots,\alpha^{k}(x_{n})])$

$=D_{\eta}D_{\xi}([x_{1},x_{2},\cdots,x_{n}])$

$-\sum^{n}_{i=2}(-1)^{(\xi+\eta)|X_{i-1}|}
[\alpha^{k+s}(x_{1}),\cdots,D_{\eta}D_{\xi}(x_{i}),\cdots,\alpha^{k+s}(x_{n})].$

Hence, it deduces that

$[[D_{\xi},D_{\eta}](x_{1}),\alpha^{k+s}(x_{2}),\cdots,\alpha^{k+s}(x_{n})]$

$=[D_{\xi}D_{\eta}(x_{1}),\alpha^{k+s}(x_{2}),\cdots,\alpha^{k+s}(x_{n})]
-(-1)^{\xi\eta}[D_{\eta}D_{\xi}(x_{1}),\alpha^{k+s}(x_{2}),\cdots,\alpha^{k+s}(x_{n})]$

$=[D_{\xi},D_{\eta}]([x_{1},x_{2},\cdots,x_{n}])$.

Similarly, we also get

$[[D_{\xi},D_{\eta}](x_{1}),\alpha^{k+s}(x_{2}),\cdots,\alpha^{k+s}(x_{n})]$

$=(-1)^{(\xi+\eta)|X_{i-1}|}
[\alpha^{k+s}(x_{1}),\cdots,[D_{\xi},D_{\eta}](x_{i}),\cdots,\alpha^{k+s}(x_{n})].$

Thus, $[D_{\xi},D_{\eta}]\in C_{\alpha^{k+s}}(N)$ of degree $\xi+\eta$. Hence, we have proven that $[Der(N),C(N)]\subseteq C(N).$

(2) Using the similar method of proving (1), it is easily obtained.

(3) Suppose that $D_{\xi}\in C_{\alpha^{k}}(N),\ D_{\eta}\in Der_{\alpha^{s}}(N)$ for any  $x_{i} \in hg(N), i\in J$. We have

$D_{\xi}D_{\eta}([x_{1},x_{2},\cdots,x_{n}])$

$=D_{\xi}(\sum^{n}_{i=1}(-1)^{\eta|X_{i-1}|}
[\alpha^{s}(x_{1}),\cdots,D_{\eta}(x_{i}),\cdots,\alpha^{s}(x_{n})])$

$=\sum^{n}_{i=1}(-1)^{(\xi+\eta)|X_{i-1}|}
[\alpha^{k+s}(x_{1}),\cdots,D_{\xi}D_{\eta}(x_{i}),\cdots,\alpha^{k+s}(x_{n})].$

Therefore, $D_{\xi}D_{\eta}\in Der_{\alpha^{k+s}}(N)$ of degree $\xi+\eta$.

Finally, we obtain
$C(N)\cdot Der (N)\subseteq  Der (N)$.

(4) Suppose that $D_{\xi}\in C_{\alpha^{k}}(N)$  for any $x_{i} \in hg(N), i\in J$, then it follows

$D_{\xi}([x_{1},x_{2},\cdots,x_{n}])$

$=(-1)^{\xi|X_{i-1}|}
[\alpha^{k}(x_{1}),\cdots,D_{\xi}(x_{i}),\cdots,\alpha^{k}(x_{n})].$

Hence, $\sum^{n}_{i=1}(-1)^{\xi|X_{i-1}|}
[\alpha^{k}(x_{1}),\cdots,D_{\xi}(x_{i}),\cdots,\alpha^{k}(x_{n})]
=nD_{\xi}([x_{1},x_{2},\cdots,x_{n}]).$

So we obtain that $D_{\xi}\in QDer_{\alpha^{k}}(N)$. Choosing
$D^{\prime}=nD_{\xi}\in C_{\alpha^{k}}(N)$, we complete the proof.

(5) Assume that $D_{\xi}\in QC_{\alpha^{k}}(N),\ D_{\eta}\in QC_{\alpha^{s}}(N)$  for any  $x_{i} \in hg(N), i\in J$. We have

$[\alpha^{k+s}(x_{1}),\cdots,[D_{\xi},D_{\eta}](x_{i}),\cdots,\alpha^{k+s}(x_{n})]$

$=(-1)^{\xi|X_{i-1}|+\eta(|X_{i-1}|-|X_{1}|)}
[D_{\xi}(\alpha^{s}(x_{1})),D_{\eta}(\alpha^{k}(x_{2})),\cdots,\alpha^{k+s}(x_{i}),\cdots,\alpha^{k+s}(x_{n})]$

$-(-1)^{\xi\eta+\eta(|X_{i-1}|-|X_{1}|)+\xi(\eta+|X_{i-1}|)}
[D_{\xi}(\alpha^{s}(x_{1})),D_{\eta}(\alpha^{k}(x_{2})),\cdots,\alpha^{k+s}(x_{i}),\cdots,\alpha^{k+s}(x_{n})]$

$=0$.

Thus, one has

$$\sum^{n}_{i=1}(-1)^{(\xi+\eta)|X_{i-1}|}
[\alpha^{k+s}(x_{1}),\cdots,[D_{\xi},D_{\eta}](x_{i}),\cdots,\alpha^{k+s}(x_{n})]=0.$$

This indicates that $[D_{\xi},D_{\eta}]\in QDer_{\alpha^{k+s}}(N)$ of degree $\xi+\eta$.

(6) It is obvious.\QED

\begin{prop}\label{po1}
We have
$QC(N)+[QC(N),QC(N)]$ is a subalgebra of $GDer(N)$.
\end{prop}

\ni\ni{\it Proof.}\ \
According to Prop~3.2(5),(6), it shows

$$QC(N)+[QC(N),QC(N)]\subseteq GDer(N)$$

and

$[QC(N)+[QC(N),QC(N)],QC(N)+[QC(N),QC(N)]]$

$\subseteq [QC(N)+QDer(N),QC(N)+[QC(N),QC(N)]]$

$\subseteq[QC(N),QC(N)]+[QC(N),[QC(N),QC(N)]]+[QDer(N),QC(N)]
+[QDer(N),[QC(N),QC(N)]]$.

Using the graded Hom-Jacobi identity, it is easy to verify that
$$[QDer(N),[QC(N),QC(N)]]\subseteq [QC(N),QC(N)],$$
 Thus,
$[QC(N)+[QC(N),QC(N)],QC(N)+[QC(N),QC(N)]]\subseteq QC(N)+[QC(N),QC(N)].$  \QED

\begin{prop}\label{po1}
Let $\alpha$ be a surjective mapping, then $[C(N),QC(N)]\subseteq Hom(N,Z(N))$. Moreover, if $Z(N)=\{0\}$, then $[C(N),QC(N)]=\{0\}$.
\end{prop}
\ni\ni{\it Proof.}\ \
 Assume that $D_{\xi}\in C_{\alpha^{k}}(N),\ D_{\eta}\in QC_{\alpha^{s}}(N)$, and $x_{1}\in hg(N)$.
Since $\alpha$ is surjective, for any  $y_{i}^{\prime}\in N$, then exists
$y_{i}\in N$ such that $y_{i}^{\prime}=\alpha^{k+s}(y_{i})$ where $2\leq i\leq n$. Thus,

$[[D_{\xi},D_{\eta}](x_{1}),y_{2}^{\prime},\cdots, y_{n}^{\prime}]$

$=[[D_{\xi},D_{\eta}](x_{1}),\alpha^{k+s}(y_{2}),\cdots, \alpha^{k+s}(y_{n})]$

$=[D_{\xi}D_{\eta}(x_{1}),\alpha^{k+s}(y_{2}),\cdots, \alpha^{k+s}(y_{n})]
-(-1)^{\xi\eta}[D_{\eta}D_{\xi}(x_{1}),\alpha^{k+s}(y_{2}),\cdots, \alpha^{k+s}(y_{n})]$

$=D_{\xi}([D_{\eta}(x_{1}),\alpha^{s}(y_{2}),\cdots, \alpha^{s}(y_{n})])
-(-1)^{\eta|X_{1}|}[D_{\xi}(\alpha^{s}(x_{1})),D_{\eta}(\alpha^{k}(y_{2})),\cdots, \alpha^{k+s}(y_{n})]$

$=D_{\xi}([D_{\eta}(x_{1}),\alpha^{s}(y_{2}),\cdots, \alpha^{s}(y_{n})])
-(-1)^{\eta|X_{1}|}D_{\xi}([\alpha^{s}(x_{1}),D_{\eta}(y_{2}),\cdots, \alpha^{s}(y_{n})])$

$=D_{\xi}([D_{\eta}(x_{1}),\alpha^{s}(y_{2}),\cdots, \alpha^{s}(y_{n})]
-(-1)^{\eta|X_{1}|}[\alpha^{s}(x_{1}),D_{\eta}(y_{2}),\cdots, \alpha^{s}(y_{n})])$

$=0$.

Hence, $[D_{\xi},D_{\eta}](x_{1})\in Z(N)$, and $[D_{\xi},D_{\eta}] \in Hom(N, Z(N))$ as desired. Furthermore,
if $Z(N)=\{0\}$, we know that $[C(N),QC(N)]=\{0\}$.   \QED

\begin{defi}\rm
An algebra $L$ over $\mathbb{F}$ is called Jordan algebra if $L$ is commutative, equipped with the following relation:~($a,u\in L$)
\begin{equation}\label{381}
a^{2}(ua)=(a^{2}u)a.
\end{equation}
The equality (3.1) implies the following equality
\begin{equation}\label{382}
((ab)u)c+((bc)u)a+((ca)u)b=(ab)(uc)+(bc)(ua)+(ca)(ub).
\end{equation}
\end{defi}
Conversely, (3.1) can be obtained from (3.2) by taking $b=c=a$ in (3.2).

\begin{defi}(\cite{AAM2}~)\rm
~A superalgebra $L=L_{\overline{0}}\oplus L_{\overline{1}}$ is called Jordan superalgebra over $\mathbb{F}$ if the multiplication satisfies:~($a,b,u,c\in hg(L)$)
$$ab=(-1)^{|a||b|}ba,$$
$$(-1)^{|c|(|a|+|u|)}((ab)u)c+(-1)^{|a|(|b|+|u|)}((bc)u)a+(-1)^{|b|(|c|+|u|)}((ca)u)b$$
\begin{equation}\label{383}
=(-1)^{|c|(|a|+|u|)}(ab)(uc)+(-1)^{|a|(|b|+|u|)}(bc)(ua)+(-1)^{|b|(|c|+|u|)}(ca)(ub).
\end{equation}
\end{defi}

The equality (3.3) is called the Super-Jordan identity.

\begin{defi}(\cite{AAM2}~)\rm
~Let $(L,\mu,\alpha)$ be a Hom-superalgebra.

(1) The Hom-associator of $L$ is a trilinear map $as_{\alpha}:L\times L\times L\rightarrow L$ defined as
$$as_{\alpha}=\mu\circ (\mu\otimes \alpha-\alpha\otimes\mu).$$
In terms of elements, the map $as_{\alpha}$ is given by:~($x,y,z\in hg(L)$)
$$as_{\alpha}(x,y,z)=\mu(\mu(x,y),\alpha(z))-\mu(\alpha(x),\mu(y,z)).$$

(2) Let $L$ be a Hom-superalgebra over  $\mathbb{F}$ with an even bilinear multiplication $\circ$. If $\alpha: L\rightarrow L$ is an even linear map, then $(L,\circ,\alpha)$ is a Hom-Jordan superalgebra if $L$ satisfies the following identities:~($x,y,z,w\in hg(L)$)
$$x\circ y=(-1)^{|x||y|}y\circ x,$$

$$(-1)^{|z|(|x|+|w|)}as_{\alpha}(x\circ y,\alpha(w),\alpha(z))
+(-1)^{|x|(|y|+|w|)}as_{\alpha}(y\circ z,\alpha(w),\alpha(x))$$
\begin{equation}\label{384}
+(-1)^{|y|(|z|+|w|)}as_{\alpha}(z\circ x,\alpha(w),\alpha(y))=0.
\end{equation}
\end{defi}

We call the equality (3.4) the Super-Hom-Jordan identity. When $\alpha=id$, then (3.4) becomes (3.3).

\begin{prop}\label{po1}
Let $(\Omega,[\cdot,\cdot],\widetilde{\alpha})$ be a Hom-Lie superalgebra over $\mathbb{F}$,
with the operation $D_{\xi}\bullet D_{\eta}=\frac{1}{2}(D_{\xi}D_{\eta}+(-1)^{\xi\eta}D_{\eta}D_{\xi})$ for $\alpha$-deriveration
$D_{\xi}, D_{\eta}\in hg(\Omega)$, then

(1) $(\Omega,\bullet, \alpha)$ is a Hom-Jordan superalgebra;

(2) $(QC(N),\bullet, \alpha)$ is also a Hom-Jordan superalgebra.
\end{prop}
\ni\ni{\it Proof.}\ \
(1) Assume that $D_{\xi},D_{\eta}\in hg(\Omega)$, we have
\begin{eqnarray*}
D_{\xi}\bullet D_{\eta}&=&\frac{1}{2}(D_{\xi}D_{\eta}+(-1)^{\xi\eta}D_{\eta}D_{\xi})\\
&=&\frac{1}{2}(-1)^{\xi\eta}(D_{\eta}D_{\xi}+(-1)^{\xi\eta}D_{\xi}D_{\eta})\\
&=&(-1)^{\xi\eta}D_{\eta}\bullet D_{\xi}.\\
\end{eqnarray*}

For one thing, it follows

$((D_{\xi}\bullet D_{\eta})\bullet \alpha(D_{\theta}))\bullet \alpha^{2}(D_{\gamma}))$

$=\frac{1}{8}\big(D_{\xi}D_{\eta}\alpha(D_{\theta})\alpha^{2}(D_{\gamma})
+(-1)^{\xi\eta}D_{\eta}D_{\xi}\alpha(D_{\theta})\alpha^{2}(D_{\gamma})
+(-1)^{(\xi+\eta)\theta}\alpha(D_{\theta})D_{\xi}D_{\eta}\alpha^{2}(D_{\gamma})$

$+(-1)^{\xi\eta+(\xi+\eta)\theta}\alpha(D_{\theta})D_{\eta}D_{\xi}\alpha^{2}(D_{\gamma})
+(-1)^{(\xi+\eta+\theta)\gamma}\alpha^{2}(D_{\gamma})D_{\xi}D_{\eta}\alpha(D_{\theta})$

$+(-1)^{\xi\eta+(\xi+\eta+\theta)\gamma}\alpha^{2}(D_{\gamma})D_{\eta}D_{\xi}\alpha(D_{\theta})
+(-1)^{(\xi+\eta)\theta+(\xi+\eta+\theta)\gamma}\alpha^{2}(D_{\gamma})\alpha(D_{\theta})D_{\xi}D_{\eta}$

$+(-1)^{\xi\eta+(\xi+\eta)\theta+(\xi+\eta+\theta)\gamma}
\alpha^{2}(D_{\gamma})\alpha(D_{\theta})D_{\eta}D_{\xi}\big)$.

For another, we have

$\alpha(D_{\xi}\bullet D_{\eta})\bullet(\alpha(D_{\theta})\bullet \alpha(D_{\gamma}))$

$=\frac{1}{8}\big(\alpha(D_{\xi}D_{\eta})\alpha(D_{\theta})\alpha(D_{\gamma})
+(-1)^{(\xi+\eta)(\theta+\gamma)}\alpha(D_{\theta})\alpha(D_{\gamma})\alpha(D_{\xi}D_{\eta})$

$+(-1)^{\theta\gamma}\alpha(D_{\xi}D_{\eta})\alpha(D_{\gamma})\alpha(D_{\theta})
+(-1)^{\theta\gamma+(\xi+\eta)(\theta+\gamma)}\alpha(D_{\gamma})\alpha(D_{\theta})\alpha(D_{\xi}D_{\eta})$

$+(-1)^{\xi\eta}\alpha(D_{\eta}D_{\xi})\alpha(D_{\theta})\alpha(D_{\gamma})
+(-1)^{\xi\eta+(\xi+\eta)(\theta+\gamma)}\alpha(D_{\theta})\alpha(D_{\gamma})\alpha(D_{\eta}D_{\xi})$

$+(-1)^{\xi\eta+\theta\gamma}\alpha(D_{\eta}D_{\xi})\alpha(D_{\gamma})\alpha(D_{\theta})
+(-1)^{\xi\eta+\theta\gamma+(\xi+\eta)(\theta+\gamma)}
\alpha(D_{\gamma})\alpha(D_{\theta})\alpha(D_{\eta}D_{\xi})\big).$

Thus, one has

$(-1)^{\gamma(\xi+\theta)}as_{\alpha}(D_{\xi}\bullet D_{\eta},\alpha(D_{\theta}),\alpha(D_{\gamma}))$

$=\frac{1}{8}(-1)^{\gamma(\xi+\theta)}\big(
(-1)^{(\xi+\eta)\theta}\alpha(D_{\theta})D_{\xi}D_{\eta}\alpha^{2}(D_{\gamma})
+(-1)^{\xi\eta+(\xi+\eta)\theta}\alpha(D_{\theta})D_{\eta}D_{\xi}\alpha^{2}(D_{\gamma})$

$+(-1)^{(\xi+\eta+\theta)\gamma}\alpha^{2}(D_{\gamma})D_{\xi}D_{\eta}\alpha(D_{\theta})
+(-1)^{\xi\eta+(\xi+\eta+\theta)\gamma}\alpha^{2}(D_{\gamma})D_{\eta}D_{\xi}\alpha(D_{\theta})$

$-(-1)^{(\xi+\eta)(\theta+\gamma)}\alpha(D_{\theta})\alpha(D_{\gamma})\alpha(D_{\xi}D_{\eta})
-(-1)^{\theta\gamma}\alpha(D_{\xi}D_{\eta})\alpha(D_{\gamma})\alpha(D_{\theta})$

$-(-1)^{\xi\eta+(\xi+\eta)(\theta+\gamma)}\alpha(D_{\theta})\alpha(D_{\gamma})\alpha(D_{\eta}D_{\xi})
-(-1)^{\xi\eta+\theta\gamma}\alpha(D_{\eta}D_{\xi})\alpha(D_{\gamma})\alpha(D_{\theta})\big)$

Therefore, we get

$(-1)^{\gamma(\xi+\theta)}as_{\alpha}(D_{\xi}\bullet D_{\eta},\alpha(D_{\theta}),\alpha(D_{\gamma}))
+(-1)^{\xi(\eta+\theta)}as_{\alpha}(D_{\eta}\bullet D_{\gamma},\alpha(D_{\theta}),\alpha(D_{\xi}))$

$+(-1)^{\eta(\gamma+\theta)}as_{\alpha}(D_{\gamma}\bullet D_{\xi},\alpha(D_{\theta}),\alpha(D_{\eta}))=0.$

(2) We need to show that $D_{\xi}\bullet D_{\eta}\in QC(N)$ for any $D_{\xi}, D_{\eta}\in hg(QC(N))$.

Assume that $x_{i}\in hg(N),i\in J$, we have

$[D_{\xi}\bullet D_{\eta}(x_{1}),\alpha^{k+s}(x_{2}),\cdots,\alpha^{k+s}(x_{n})]$

$=\frac{1}{2}[D_{\xi}D_{\eta}(x_{1}),\alpha^{k+s}(x_{2}),\cdots,\alpha^{k+s}(x_{n})]
+\frac{1}{2}(-1)^{\xi\eta}[D_{\eta}D_{\xi}(x_{1}),\alpha^{k+s}(x_{2}),\cdots,\alpha^{k+s}(x_{n})]$

$=\frac{1}{2}(-1)^{\xi\eta+(\xi+\eta)|X_{i-1}|}[\alpha^{k+s}(x_{1}),\cdots,
D_{\eta}D_{\xi}(x_{i}),\cdots,\alpha^{k+s}(x_{n})]$

$+\frac{1}{2}(-1)^{(\xi+\eta)|X_{i-1}|}[\alpha^{k+s}(x_{1}),\cdots,
D_{\xi}D_{\eta}(x_{i}),\cdots,\alpha^{k+s}(x_{n})]$

$=(-1)^{(\xi+\eta)|X_{i-1}|}[\alpha^{k+s}(x_{1}),\cdots,
D_{\xi}\bullet D_{\eta}(x_{i}),\cdots,\alpha^{k+s}(x_{n})]$.

Thus, $D_{\xi}\bullet D_{\eta}\in QC(N)$.

So we have completed the proof.  \QED

\begin{prop}\label{po1}
Let $(N,[\cdot,\cdots,\cdot],\alpha)$ be a multiplicative n-Hom Lie superalgebra. Then the following statements hold:

(1) $QC(N)$ is a Lie superalgebra with $[D_{\xi},D_{\eta}]=D_{\xi}D_{\eta}-(-1)^{\xi\eta}D_{\eta}D_{\xi}$ if and only if $QC(N)$ is also a Hom-associative superalgebra;

(2) If $Z(N)=\{0\}$, then $QC(N)$ is a Lie superalgebra if and only if $[QC(N),QC(N)]=\{0\}$.
\end{prop}

\ni\ni{\it Proof.}\ \
(1)$(\Rightarrow)$ Note that $D_{\xi}D_{\eta}=D_{\xi}\bullet D_{\eta}+\frac{[D_{\xi},D_{\eta}]}{2}$. By Prop~3.8, then $D_{\xi}\bullet D_{\eta}\in QC(N),[D_{\xi},D_{\eta}]\in QC(N)$. It follows $D_{\xi}D_{\eta}\in QC(N)$.

 $(\Leftarrow) $ For any $D_{\xi}\in QC_{\alpha^{k}}(N),\ D_{\eta}\in QC_{\alpha^{s}}(N)$, then we have $D_{\xi}D_{\eta}\in QC_{\alpha^{k+s}}(N)$ and $D_{\eta}D_{\xi}\in QC_{\alpha^{k+s}}(N)$. We also see that $[D_{\xi},D_{\eta}]=D_{\xi}D_{\eta}-(-1)^{\xi\eta}D_{\eta}D_{\xi}\in QC_{\alpha^{k+s}}(N)$. Hence, $QC(N)$ is a Lie superalgebra.

(2) $(\Rightarrow)$ Let $D_{\xi}\in QC_{\alpha^{k}}(N),\ D_{\eta}\in QC_{\alpha^{s}}(N)$. Since $QC(N)$ is a Lie superalgebra, then we have

$[[D_{\xi},D_{\eta}](x_{1}),\alpha^{k+s}(x_{2}),\cdots,\alpha^{k+s}(x_{n})]$

$=(-1)^{(\xi+\eta)|X_{i-1}|}
[\alpha^{k+s}(x_{1}),\cdots,[D_{\xi},D_{\eta}](x_{i}),\cdots,\alpha^{k+s}(x_{n})].$

From the proof of Prop 3.2(5), then we obtain

$$\sum_{i=1}^{n}(-1)^{(\xi+\eta)|X_{i-1}|}
[\alpha^{k+s}(x_{1}),\cdots,[D_{\xi},D_{\eta}](x_{i}),\cdots,\alpha^{k+s}(x_{n})]=0.$$
Hence, $n[[D_{\xi},D_{\eta}](x_{1}),\alpha^{k+s}(x_{2}),\cdots,\alpha^{k+s}(x_{n})]=0$.

Since $char~\mathbb{F}=0$, we infer that
$[[D_{\xi},D_{\eta}](x_{1}),\alpha^{k+s}(x_{2}),\cdots,\alpha^{k+s}(x_{n})]=0$, i.e. $[D_{\xi},D_{\eta}]=0$.

$(\Leftarrow) $ It is clear.
  \QED

\section{Quasiderivations of  multiplicative n-Hom Lie superalgebras}

In this section, we shall investigate the quasiderivations of the multiplicative n-Hom Lie superalgebra $N$. We obtain that
$QDer(N)$ can be embedded as derivations in a larger multiplicative n-Hom Lie superalgebra. Moreover, we get that $Der(\breve{N})$ has a direct sum decomposition  when  $Z(N)=\{0\}$.

\begin{prop}\label{po1}
Let $(N,[\cdot,\cdots,\cdot],\alpha)$ be a multiplicative n-Hom Lie superalgebra over $\mathbb{F}$ and $t$ be an indeterminate. We define $\breve{N}:=\{\sum(x\otimes t+y\otimes t^{n})\ |\ x,y\in hg(N)\}, \breve{\alpha}(\breve{N}):=\{\sum(\alpha(x)\otimes t+\alpha(y)\otimes t^{n})\ |\ x,y\in hg(N)\}$. Then $\breve{N}$ is a multiplicative n-Hom Lie superalgebra with the bracket:
$$[x_{1}\otimes t^{j_{1}},x_{2}\otimes t^{j_{2}},\cdots,x_{n}\otimes t^{j_{n}}]
=[x_{1},x_{2},\cdots,x_{n}]\otimes t^{\sum_{k=1}^{n} j_{k}},$$
for any $x_{i}\in hg(N), j_{k}\in \{1,n\},i,k\in J$.
\end{prop}

\ni\ni{\it Proof.}\ \
For any $x_{i},y_{i}\in hg(N)$ where $i\in J$, then one has

$[x_{1}\otimes t^{j_{1}},\cdots,x_{i}\otimes t^{j_{i}},x_{i+1}\otimes t^{j_{i+1}},\cdots,x_{n}\otimes t^{j_{n}}]$

$=[x_{1},\cdots,x_{i},x_{i+1},\cdots,x_{n}]\otimes t^{\sum_{k=1}^{n} j_{k}}$

$=-(-1)^{|x_{i}||x_{i+1}|}[x_{1},\cdots,x_{i+1},x_{i},\cdots,x_{n}]\otimes t^{\sum_{k=1}^{n} j_{k}}$

$=-(-1)^{|x_{i}||x_{i+1}|}[x_{1}\otimes t^{j_{1}},\cdots,x_{i+1}\otimes t^{j_{i+1}},x_{i}\otimes t^{j_{i}},\cdots,x_{n}\otimes t^{j_{n}}].$

It shows that

$[\breve{\alpha}(x_{1}\otimes t^{j_{1}}),\cdots,\breve{\alpha}(x_{n-1}\otimes t^{j_{n-1}}),[y_{1}\otimes t^{l_{1}},\cdots,y_{n}\otimes t^{l_{n}}]]$

$=\sum^{n}_{i=1}(-1)^{|X_{n-1}||Y_{i-1}|}
[\alpha(y_{1}),\cdots,\alpha(y_{i-1}),[x_{1},\cdots,x_{n-1},y_{i}],
\alpha(y_{i+1}),\cdots,\alpha(y_{n})]\otimes t^{(\sum_{k=1}^{n-1}j_{k}+\sum_{k=1}^{n}l_{k})}$

$=\sum^{n}_{i=1}(-1)^{|X_{n-1}||Y_{i-1}|}
[\breve{\alpha}(y_{1}\otimes t^{l_{1}}),\cdots,[x_{1}\otimes t^{j_{1}},\cdots,x_{n-1}\otimes t^{j_{n-1}},y_{i}\otimes t^{l_{i}}],\cdots,\breve{\alpha}(y_{n}\otimes t^{l_{n}})].$

Hence, $\breve{N}$ is a multiplicative n-Hom Lie superalgebra.  \QED

In the following we shall investigate some properties of generalized derivations of $\breve{N}$. For convenience, we denote $xt(xt^{n})$ by $x\otimes t(x\otimes t^{n})$.
If $U$ is a $\mathbb{Z}_{2}$-graded subspace of $N$ such that $N=U\oplus [N,\cdots,N]$, then
$$\breve{N}=Nt+Nt^{n}=Nt+Ut^{n}+[N,\cdots,N]t^{n}.$$

Now we define a map $\varphi: QDer(N)\rightarrow End(\breve{N})$ satisfying
$$\varphi(D)(at+ut^{n}+bt^{n})=D(a)t+D^{\prime}(b)t^{n},$$
where $D\in QDer(N)$, $a\in hg(N), u\in hg(U), b\in hg([N,\cdots,N])$ and $|a|=|b|=|u|$. In addition, $D^{\prime}$ is a linear map related to $D$ in Def~2.8.

\begin{prop}\label{po1}
Let $N,\breve{N},\varphi$ be as above. Then

(1) $|\varphi|=0$;

(2) $\varphi$ is injective and $\varphi(D)$ does not depend on the choice of $D^{\prime}$;

(3) $\varphi(QDer(N))\subseteq Der(\breve{N})$.
\end{prop}
\ni\ni{\it Proof.}\ \
(1) In view of the definition of $\varphi$, there is nothing to prove.

(2) If $\varphi(D_{1})=\varphi(D_{2})$, then it follows that

$$\varphi(D_{1})(at+ut^{n}+bt^{n})=\varphi(D_{2})(at+ut^{n}+bt^{n}).$$
In other words,
$$D_{1}(a)t+D^{\prime}_{1}(b)t^{n}=D_{2}(a)t+D^{\prime}_{2}(b)t^{n}.$$

Thus, $D_{1}(a)=D_{2}(a)$. So we obtain that $\varphi$ is injective.

Suppose there exists $D^{\prime\prime}\in End(N)$ such that
$$\varphi(D)(at+ut^{n}+bt^{n})=D(a)t+D^{\prime\prime}(b)t^{n},$$
and

$$\sum^{n}_{i=1}(-1)^{|D| |X_{i-1}|}
[\alpha^{k}(x_{1}),\cdots,\alpha^{k}(x_{i-1}),D(x_{i}),\alpha^{k}(x_{i+1}),\cdots,\alpha^{k}(x_{n})]
=D^{\prime\prime}([x_{1},\cdots,x_{n}]),$$
then we have
$$D^{\prime}([x_{1},\cdots,x_{n}])=D^{\prime\prime}([x_{1},\cdots,x_{n}]).$$

Therefore, $D^{\prime}=D^{\prime\prime}$,
%. %Hence
%$$\varphi(D)(at+ut^{n}+bt^{n})=D(a)t+D^{\prime}(b)t^{n}=D(a)t+D^{\prime\prime}(b)t^{n},$$
which implies $\varphi(D)$ is uniquely determined by $D$.

(3) We know that $[x_{1}t^{j_{1}},\cdots,x_{n}t^{j_{n}}]=[x_{1},\cdots,x_{n}]t^{\sum_{k=1}^{n} j_{k}}=0$ for any $\sum_{k=1}^{n} j_{k}\geq n+1$. Thus, we only need to check the following equation:
$$\varphi(D)([x_{1}t,\cdots,x_{n}t])=\sum^{n}_{i=1}(-1)^{|D| |X_{i-1}|}[\breve{\alpha}^{k}(x_{1}t),\cdots,\varphi(D)(x_{i}t),
\cdots,\breve{\alpha}^{k}(x_{n}t)].$$

For any $x_{i}\in hg(N),i\in J$, we have

$\varphi(D)([x_{1}t,\cdots,x_{i}t,\cdots,x_{n}t])$

$=\varphi(D)([x_{1},\cdots,x_{i},\cdots,x_{n}]t^{n})$

$=D^{\prime}([x_{1},\cdots,x_{n}])t^{n}$

$=\sum^{n}_{i=1}(-1)^{|D| |X_{i-1}|}
[\alpha^{k}(x_{1}),\cdots,\alpha^{k}(x_{i-1}),D(x_{i}),\alpha^{k}(x_{i+1}),\cdots,\alpha^{k}(x_{n})]t^{n}$

$=\sum^{n}_{i=1}(-1)^{|D| |X_{i-1}|}
[\alpha^{k}(x_{1})t,\cdots,\alpha^{k}(x_{i-1})t,D(x_{i})t,\alpha^{k}(x_{i+1})t,\cdots,\alpha^{k}(x_{n})t]$

$=\sum^{n}_{i=1}(-1)^{|D| |X_{i-1}|}
[\breve{\alpha}^{k}(x_{1}t),\cdots,\breve{\alpha}^{k}(x_{i-1}t),\varphi(D)(x_{i}t),
\breve{\alpha}^{k}(x_{i+1}t),\cdots,\breve{\alpha}^{k}(x_{n}t)].$

Hence, we get $\varphi(D)\in Der(\breve{N})$ for any $D\in QDer(N)$.   \QED

\begin{prop}\label{po1}
Let $N$ be a multiplicative n-Hom Lie superalgebra such that $Z(N)=\{0\}$ and $\breve{N},\varphi$ be as defined above. Then
$$Der(\breve{N})=\varphi(QDer(N))\oplus ZDer(\breve{N}).$$
\end{prop}
\ni\ni{\it Proof.}\ \
Due to $Z(N)=\{0\}$, then one has $Z(\breve{N})=Nt^{n}$. For any $g \in Der(\breve{N})$, we easily see that $g(Z(\breve{N}))\subseteq Z(\breve{N})$. Thus, $g(Ut^{n})\subseteq g(Z(\breve{N}))\subseteq Z(\breve{N})=Nt^{n}$. Now define a map $f:Nt+Ut^{n}+[N,\cdots,N]t^{n}\rightarrow Nt^{n}$ by

\begin{equation*}\label{111111}
f(x)=\left\{\begin{array}{lll}
g(x)\cap Nt^n&\mbox{if \ }x\in Nt,\\
g(x)&\mbox{if \ }x\in Ut^n,\\
0&\mbox{if \ }x\in [N,\cdots,N]t^{n}.
\end{array}\right.
\end{equation*}
It is is easy to see that $f$ is a linear map. Note that
$$f([\breve{N},\cdots,\breve{N}])=f([N,\cdots,N]t^{n})=0,$$
$$[f(\breve{N}),\breve{\alpha}^{k}(\breve{N}),\cdots,\breve{\alpha}^{k}(\breve{N})]\subseteq [Nt^{n},\alpha^{k}(N)t+\alpha^{k}(N)t^{n},\cdots,\alpha^{k}(N)t+\alpha^{k}(N)t^{n}]=0,$$
hence, $f\in ZDer(\breve{N})$.

Now set $h=g-f$, then

$$h(Nt)=g(Nt)-g(Nt)\cap Nt^{n}=g(Nt)-Nt^{n}\subseteq Nt,\ h(Ut^{n})=0.$$

We also obtain

$$h([N,\cdots,N]t^{n})=g([\breve{N},\cdots,\breve{N}])\subseteq [\breve{N},\cdots,\breve{N}]=[N,\cdots,N]t^{n},$$
there exist $D,D^{\prime}\in End(N)$ such that for any $a\in hg(N), b\in hg([N,\cdots,N])$,

$$h(at)=D(a)t,\ h(bt^{n})=D^{\prime}(b)t^{n}.$$

Since $h\in Der(\breve{N})$ and from the definition of $Der(\breve{N})$, then we get
$$\sum^{n}_{i=1}(-1)^{|h| |A_{i-1}|}
[\breve{\alpha}^{k}(a_{1}t),\cdots,\breve{\alpha}^{k}(a_{i-1}t),h(a_{i}t),
\breve{\alpha}^{k}(a_{i+1}t),\cdots,\breve{\alpha}^{k}(a_{n}t)]
=h([a_{1}t,\cdots,a_{n}t]),$$
for any $a_{i}\in hg(N),i\in J$. Hence,
$$\sum^{n}_{i=1}(-1)^{|D| |A_{i-1}|}[\breve{\alpha}^{k}(a_{1}),\cdots,D(a_{i}),\cdots,
\breve{\alpha}^{k}(a_{n})]=D^{\prime}([a_{1},\cdots,a_{n}]).$$
Thus, $D\in QDer(N)$. Furthermore, $h=\varphi(D)\in \varphi(QDer(N))$. Thus, $Der(\breve{N})\subseteq \varphi(QDer(N))+ ZDer(\breve{N}).$ By Prop 4.2(3), we have $Der(\breve{N})=\varphi(QDer(N))+ ZDer(\breve{N}).$

For any $f\in \varphi(QDer(N))\cap ZDer(\breve{N})$, then there exists $D\in QDer(N)$ such that $f=\varphi(D)$. Then,
$$f(at+ut^{n}+bt^{n})=\varphi(D)(at+ut^{n}+bt^{n})=D(a)t+D^{\prime}(b)t^{n}$$
for any $a\in hg(N), b\in hg([N,\cdots,N])$.

We also know that $f\in ZDer(\breve{N})$, then one shows
$$f(at+ut^{n}+bt^{n})\in Z(\breve{N})=Nt^{n}.$$

This means that $D(a)=0$ for any $a\in hg(N)$. Furthermore, we get $D=0$. Thus, $f=\varphi(0)=0$.

It follows that $Der(\breve{N})=\varphi(QDer(N))\oplus ZDer(\breve{N})$.    \QED

\end{document}